\documentstyle[amstex,amscd,12pt]{article}
\makeatletter
\date{}
\newtheorem{theorem}{Theorem}[section]

\newtheorem{proposition}[theorem]{Proposition}

\newcommand{\edim}{{\rm e-dim}}
\newcommand{\z}{{\Bbb Z}}

\newcommand{\s}{{\Bbb S}}
\newcommand {\h}{{\check H}}

 \newcommand{\invlim}{\raisebox{-1ex}{$\stackrel{\hbox{lim}}{\leftarrow}$}}

\newcommand{\lo}{\longrightarrow}

\newcommand{\tor}{{\rm Tor}}
\begin{document}

\title{Universal acyclic  resolutions for finitely generated coefficient groups\\
}

\author{Michael  Levin}

\maketitle
\begin{abstract}
We prove that for every compactum $X$   and every integer $n \geq 2$
there  are  a compactum $Z$ of $\dim \leq n$ and a surjective $UV^{n-1}$-map
$r: Z \lo X$   having the  property  that:

 for every finitely generated abelian group $G$   and every integer $k \geq 2$
  such that  $\dim_G  X \leq k \leq n$ we have      $\dim_G Z \leq k$ and
  $r$ is $G$-acyclic,
 or equivalently:

  for every simply connected
 CW-complex $K$  with finitely generated homotopy groups    such that
   $\edim X \leq K$ we have $\edim  Z \leq K$ and $r$ is $K$-acyclic.
(A space is $K$-acyclic
    if every map from the space to $K$ is null-homotopic.  A map is $K$-acyclic
   if  every fiber is $K$-acyclic.)
\bigskip
\\
{\bf Math. Subj. Class.:} 55M10, 54F45.
\end{abstract}
\begin{section}{Introduction}
A space $X$  is always assumed to be separable metrizable. The cohomological dimension $\dim_G X$
of  $X$  with respect to an abelian group $G$ is   the least number $n$ such that $\h^{n+1}
(X,A;G)=0$ for every closed subset $A$ of $X$.
 The covering dimension
 $\dim X $ coincides with the integral cohomological dimension $\dim_\z X$
 if $X$ is finite dimensional.  In 1987 Dranishnikov \cite{dr0} showed that     there is
an infinite dimensional compactum (=compact metric space)
 with
 finite integral cohomological dimension.  In 1978   Edwards discovered his
 resolution theorem  \cite{ed1, w1}
which shows  that   despite Dranishnikov's  example
 the dimension functions $\dim$ and $\dim_\z$ are closely related even for
 infinite dimensional compacta. The    Edwards  resolution
 theorem says that a compactum of $\dim_\z \leq n$ can be obtained
 as the  image of  a cell-like map defined  on a compactum of  $\dim \leq n$.
   A  compactum    $X$ is cell-like if any map $f : X \lo K$ from $X$ to a CW-complex  $K$  is
   null-homotopic.
   A  map is cell-like if its fibers are cell-like. The reduced $\rm {\check C}$ech
 cohomology groups of  a cell-like compactum   are trivial with respect to any group $G$.

 Below are stated a few theorems  representing main directions of generalizing
the   Edwards  resolution  theorem.  We recall  some notions used in the theorems.

A space is    $G$-acyclic if its       reduced $\rm {\check C}$ech   cohomology groups modulo $G$
  are trivial,   a map is $G$-acyclic if
   every  fiber is  $G$-acyclic.
     By the Vietoris-Begle  theorem a $G$-acyclic map  of compacta
cannot raise the cohomological dimension $\dim_G$.

  A compactum $X$ is  approximately $n$-connected  if any embedding of $X$ into an ANR
  has the $UV^n$-property,   i.e.  for every    neighborhood $ U$ of $X$ there is
  a smaller neighborhood  $X \subset V \subset U$  such that the inclusion $V \subset U$
   induces
  the zero homomorphism of the homotopy groups in $\dim \leq n$.
  An     approximately $n$-connected  compactum has trivial    reduced $\rm {\check C}$ech
 cohomology groups     in $\dim \leq n$   with respect to any group $G$.
A map is called a    $UV^n$-map if  every fiber is      approximately $n$-connected.

 Let $K$ be a CW-complex.       A space is said to be $K$-acyclic if every map from the space
 to $K$ is null-homotopic and
 a map   is said to be $K$-acyclic
 if every fiber  is  $K$-acyclic.

  \begin{theorem} {\rm (\cite{dr1})}
\label{dr1}
Let $p$ be  a prime number and
 let $X$ be  a compactum with $\dim_{\z_p} X \leq n$.  Then there are a compactum
 $Z $  with
  $\dim Z \leq n$ and a $\z_p$-acyclic $UV^{n-1}$-map $r : Z \lo X$  from $Z$ onto $X$.
\end{theorem}

      \begin{theorem}  {\rm (\cite{l1})}
\label{l1}
 Let $G$ be an abelian group and
 let $X$ be  a compactum with $\dim_G X \leq n$, $ n \geq 2$.  Then there are a compactum
 $Z $  with
 $\dim_G Z \leq n$ and  $\dim Z \leq n+1$ and a $G$-acyclic map $r : Z \lo X$  from $Z$ onto $X$.
\end{theorem}

   \begin{theorem}   {\rm(\cite{l3})}
 \label{l3}
 Let $X$ be a compactum  with $\dim_\z X \leq n \geq 2$.  Then there exist
 a compactum $Z$ with $\dim Z \leq n$  and a cell-like map $r : Z \lo X$ from
 $Z$ onto $X$ such that for every integer $k \geq 2$ and every group $G$  such that
 $\dim_G X \leq k$ we have $\dim_G Z \leq k$.
 \end{theorem}

   Each of these theorems emphasizes   different aspects of the Edwards resolution.
  The special feature of    Theorem \ref{l3}  is  that  it     takes care of all possible cohomological
  dimensions   bigger than $1$.   It is very natural   to try to generalize  Theorems \ref{dr1}
  and \ref{l1} in the spirit   of       Theorem \ref{l3}.
  As a part of this project this paper is mainly  devoted   to proving
    the following version of  Theorem \ref{dr1}.

\begin{theorem}
\label{t1}
Let $X$ be a compactum.  Then for every integer   $n \geq 2$
  there  are  a compactum $Z$ of $\dim \leq n$ and a surjective $UV^{n-1}$-map
$r: Z \lo X$   having the  property  that
 for every finitely generated abelian group $G$   and every integer $k \geq 2$
  such that  $\dim_G  X \leq k \leq n$ we have  that    $\dim_G Z \leq k$ and
  $r$ is $G$-acyclic.
  \end{theorem}

 Theorem \ref{t1}  can be reformulated in terms of extensional dimension \cite{drdyd0, drdyd1}.
 The extensional dimension of $X$     is said
  not to exceed  a CW-complex $K$, written $\edim X \leq K$,  if for every closed subset $A$
 of $X$ and every map $f : A \lo K$ there is an extension of $f$ over $X$.
 It is well-known that $\dim X \leq n$ is equivalent to $\edim X \leq \s^n$
 and $\dim_G X \leq  n$ is equivalent to $\edim X \leq K(G,n)$ where
 $K(G,n)$ is an Eilenberg-Mac Lane complex of type $(G,n)$.

 The following theorem  shows a close connection between  cohomological
 and extensional dimensions.
\begin{theorem}
%\label{dr3}
{\rm (\cite{dr3})}
\label{dr3}
Let $X$ be a  compactum and let
$K$ be a simply connected CW-complex. Consider the following
conditions:

(1) $\edim X \leq  K$;

(2) $\dim_{H_i(K)} X \leq i$ for every $i > 1$;

(3) $\dim_{\pi_i(K)} X \leq i$ for every $i >1$.

Then   (2) and (3) are equivalent and (1) implies both (2) and (3).
If $X$ is  finite dimensional then all the conditions are equivalent.
\end{theorem}

Theorems \ref{t1} and \ref{dr3} imply
\begin{theorem}
 \label{t2}
 Let $X$ be a compactum.    Then  for every integer $n\geq 2$ there exist
 a compactum $Z$ with $\dim Z \leq n$  and a surjective $UV^{n-1}$-map $r : Z \lo X$
  such that  for every  simply connected CW-complex $K$   with finitely generated
  homotopy groups
   such that
 $\edim  X \leq K$ we have $\edim Z \leq K$  and $r$ is $K$-acyclic.
 \end{theorem}
 {\bf Proof.} Let $Z$ and $r : Z \lo X$ be as in Theorem \ref{t1}.
 Let a simply connected CW-complex $K$ with finitely generated homotopy groups   be such that
 $\edim X \leq K$. Then ${H_*(K)}$
 are finitely generated and by Theorem \ref{dr3},  $\dim_{H_i(K)} X \leq i$
 and  $\dim_{\pi_i(K)} X \leq i$ for every $i > 1$.
 Hence by Theorem  \ref{t1}, $\dim_{H_i(K)} Z \leq i$ for every $1 < i \leq n$.
 Since $\dim Z \leq n$,  $\dim_{H_i(K)} Z \leq i$ for $i > n$ and
 it follows from   Theorem \ref{dr3} that
 $\edim Z \leq K$.

  Let $x \in X$ and  let $Cr^{-1}(x)$  and $\Sigma r^{-1}(x)$
 be the cone and the suspension of $r^{-1}(x)$   respectively.
 Then  $\h^{i+1}(Cr^{-1}(x),r^{-1}(x);\pi_i{(K)})= \h^{i+1}(Cr^{-1}(x)/r^{-1}(x);\pi_i{(K)})=
 \h^{i+1}(\Sigma r^{-1}(x);\pi_i{(K)})
 =\Tilde{\h}^{i}(r^{-1}(x);\pi_i{(K}))$.  Recall that    $\dim_{\pi_i(K)} X \leq i$ for every $i > 1$ and
 hence,   by Theorem \ref{t1}, $r$ is     $\pi_i(K)$-acyclic  for $ 1< i \leq n $.
 Then, since $\dim Z \leq  n$,   $\h^{i}(r^{-1}(x);\pi_i{(K)})=0$
     for  every $ i  > 1 $ and therefore
 ${\h}^{i+1}(Cr^{-1}(x),r^{-1}(x);\pi_i{(K)})=0$   for $i > 1$.
 Since    $Cr^{-1}(x)$ is finite dimensional     it follows  from
  Obstruction Theory      that every map from    $r^{-1}(x)$ to $K$
  extends over $ Cr^{-1}(x)$ and hence is null-homotopic.   Thus we showed that $r$ is
  $K$-acyclic.
  \hfill $\Box$
   \\

  Theorem \ref{t2} is equivalent to Theorem  \ref{t1}. To show this we need to show
  that    Theorem \ref{t2}  implies     Theorem  \ref{t1}. Let $X$ be a compactum, let
  $n \geq 2$ and    let $r:  Z \lo X$   satisfy the conclusions  of      Theorem \ref{t2}.
  Then for every finitely generated
  abelian group $G$ with  $\dim_G  X  \leq k \leq n$, $k \geq 2$ we have
  that $\dim_G Z \leq k$.
 Let $x \in X$. The  $K(G,k)$-acyclicity of $r$    implies that $\h^{k}(r^{-1}(x); G)=0$,
 the fact    that $r$ is $UV^{n-1}$ implies that $  \Tilde{\h}^{i}(r^{-1}(x); G)=0$
 for $i \leq k-1$ and finally  $\dim_G Z \leq k$   implies that   $\h^{i}(r^{-1}(x); G)=0$ for $ i > k$.
 All this together implies that $r^{-1}(x)$ is $G$-acyclic  and hence
  Theorems \ref{t2} and  \ref{t1}       are equivalent.
\\

     It is known that      Theorem \ref{t1} does not hold for arbitrary  groups $G$
      even  if we do not require that  $r$ is   $UV^{n-1}$,  see \cite{ky2}.
    However allowing  in Theorem \ref{t1}  for the dimension of $Z$ to  be  raised to $ n+1$
     one can  drop the restrictions on $G$
      and obtain the following  theorem which we will state without a proof.

   \begin{theorem}
\label{t3}
Let $X$ be a compactum.  Then for every integer   $n \geq 2$
  there  are  a compactum $Z$ of $\dim \leq n+1$ and a surjective $UV^{n-1}$-map
$r: Z \lo X$   having the  property  that
 for every  abelian group $G$   and every integer $k \geq 2$
  such that  $\dim_G  X \leq k \leq n$ we have  that    $\dim_G Z \leq k$ and
  $r$ is $G$-acyclic.
  \end{theorem}

   Finally let us note that it would be interesting to know if the restriction
   $k \geq 2$ in Theorem \ref{t1} can be omitted.

 \end{section}
 \begin{section}{Preliminaries}
     A  map between CW-complexes  is said to be  combinatorial if  the preimage of
   every subcomplex
   of the range is a subcomplex of the domain.

  Let $M$ be  a simplicial complex and let  $M^{[k]}$        be
    the $k$-skeleton of $M$ (=the union of all simplexes of $M$ of $\dim \leq k$).
By
 a resolution $EW(M,k)$   of $M$   we mean a CW-complex $EW(M,k)$ and
 a combinatorial map
 $\omega : EW(M,k) \lo  M$ such that $\omega$ is 1-to-1 over $M^{[k]}$.
 Let $f : N \lo K$  be a map of a subcomplex $N$ of $M$ into a CW-complex $K$.
The resolution is said to be suitable for $f$   if
 the map  $f \circ\omega|_{\omega^{-1}(N)}$ extends
   to a map  $ f': EW(M,k) \lo K$.  We call $f'$  a resolving map for $f$.
     The resolution is said to be
   suitable  for  a compactum $X$
if for   every simplex $\Delta$ of $M$,
 $\edim X  \leq \omega^{-1}(\Delta)$.
     Note that if $\omega: EW(M,k) \lo M$ is a resolution suitable
 for $X$   then  for every map  $\phi :  X \lo  M$  there is  a map   $\psi : X \lo EW(M,k)$
 such that  for every simplex $\Delta$ of $M$,
  $(\omega \circ \psi)(\phi^{-1}(\Delta)) \subset \Delta$.
 We  call $\psi$ a combinatorial lifting of $\phi$.

Let $M$   be a finite simplicial complex and
 let       $f : N  \lo K$ be a cellular  map  from a subcomplex $N$
of $M$  to a CW-complex    $K$ such that
$M^{[k]}\subset N$.
A standard way of constructing
a  resolution  suitable for $f$ is described in \cite{l3}.
Such a  resolution  $\omega: EW(M,k) \lo M$  is called the standard resolution    of $M$
for $f$ and it has the following properties:

$\omega$ is a map onto and for every simplex $\Delta$ of $M$, $\omega^{-1}(\Delta)$ is  either
contractible  or homotopy equivalent to $K$;

the (integral) homology groups  of        $EW(M,k)$
are finitely generated if so are the homology groups of $K$.  This property can be derived from
the previous one   using  the Mayer-Vietoris sequence and
induction on the number of simplexes of $M$;

 $EW(M,k)$ is $(k-1)$-connected if so are $M$ and $K$;

for every subcomplex $T$ of $M$,  $\omega|_{\omega^{-1}(T)} : EW(T,k)= \omega^{-1}(T) \lo T$
is the standard resolution of $T$ for $f|_{N \cap T} :   N \cap T \lo K$.
     \\

 All  groups are assumed to be
   abelian  and functions between groups are homomorphisms.
   $\cal P$ stands for the set of  primes.       Let $G$ be a group  and $p \in \cal P$.
   We say that $g \in G$ is $p$-torsion if  $p^k g =0$ for some  integer $k \geq 1$.
$\tor_p G$ stands for the subgroup of the $p$-torsion elements of $G$.
   $G$ is $p$-torsion if      $\tor_p G=G$,  $G$ is  $p$-torsion free if    $\tor_p G=0$  and
   $G$ is $p$-divisible if  for
   every  $g \in G$  there is $h\in G$ such that $ph=g$.
      $G$ is $p$-local if it is $q$-divisible and $q$-torsion free
   for every $q \in {\cal P} $, $q \neq p$.\\

    Let $G$ be a group,  let $\alpha : L \lo M$ be a surjective  combinatorial map of a CW-complex
$L$ and a finite simplicial complex $M$ and let $n$ be a positive integer
 such that  $\Tilde{H}_i(\alpha^{-1}(\Delta) ; G)=0$
for every $i < n$ and  every  simplex $\Delta $ of $M$.     One can show by induction on
the number of simplexes of $M$  using  the Mayer-Vietoris sequence and the Five Lemma
that $\alpha_* :   \Tilde{H}_i (L; G) \lo    \Tilde{H}_i (M; G)$ is an isomorphism
for $i < n$.   We will refer to this fact as the combinatorial Vietoris-Begle theorem.  \\

   We need  the following slightly more precise version   of Proposition 2.1, (i) of \cite{l3}.

  \begin{proposition}
  \label{p1}
  Let $2 \leq k \leq n $, $p \in \cal P$  and let $M$ be an $(n-1)$-connected finite simplicial complex.
  Let $\omega : EW(M,k)\lo M$ be the standard resolution  for
 a cellular map  $f : N \lo K(\z_p,k)$  from   a subcomplex $N$  of $M$ containing $M^{[k]}$.
  Then
  $\pi_i (  EW(M,k))$ is
  $p$-torsion    for every $1\leq i \leq n-1$.
\end{proposition}
  {\bf Proof.}
  Recall that   $\omega$ is a combinatorial surjective map and
    for every  simplex $\Delta$ of $ M$,
  $\omega^{-1}(\Delta)$ is either contractible  or homotopy equivalent to  $K(\z_p,k)$.
 $EW(M,k)$ is $(k-1)$-connected since so are $M$ and        $K(\z_p,k)$, and
  ${H}_i( EW(M,k))$  are finitely generated since  so are  ${H}_i( K(\z_p,k))$.

    By the generalized Hurewicz theorem the groups     $H_i ( K(\z_p, k))$, $i\geq 1$
  are $p$-torsion and therefore    $H_i ( K(\z_p, k))$, $i\geq 1$ are $p$-local.
   Let $q \in \cal  P$  and $q \neq p$.
 From  the $p$-locality of   $H_i ( K(\z_p, k))$, $i\geq 1$    it follows  that
          $H_i ( K(\z_p , k);\z_q)=0$, $i\geq 1$.
   Since  $M$ is $(n-1)$-connected,
   the combinatorial Vietoris-Begle   theorem    implies that
    $   H_i(EW(M,k); \z_q)=0$,  $1\leq i \leq n-1$.
 Then from the universal coefficient theorem it follows that
 $H_i(EW(M,k))\otimes \z_q=0$  for $1\leq i \leq n-1$  and every $q \in \cal P$, $q \neq p$.
Since  $H_i(EW(M,k))$ is  finitely generated,
 the  last property implies  that      $  H_i(EW(M,k))$,  $1\leq i \leq n-1$ is
 $p$-torsion   and by the generalized Hurewicz theorem
  $  \pi_i (EW(M,k))$,  $1\leq i \leq  n-1$   is       $p$-torsion.   \\
  \hfill $\Box$
     \\\\
 In the proof of Theorem \ref{t1}    we will also  use the   following facts.
  \begin{proposition}  {\rm (\cite{l3})}
 \label{p2}
 Let $K$ be a simply connected CW-complex such that  $K$ has only finitely many non-trivial
 homotopy groups.   Let $X$ be a compactum such that $\dim_{\pi_i (K)} X  \leq i$
 for  $i>1$.  Then $\edim X \leq K$.
 \end{proposition}

  Let $K'$ be a simplicial complex.  We say that    maps
  $h : K \lo K'$, $g : L \lo L'$,  $\alpha : L \lo K$ and $\alpha' : L' \lo K'$
 \[
    \begin{CD}
           L @>\alpha>>   K\\
           @V g VV        @V h VV \\
             L' @>\alpha'>>   K'
    \end{CD}
\]
 combinatorially commute   if   for every
 simplex  $\Delta$ of $K'$ we have that
 $(\alpha' \circ g)(( h\circ \alpha)^{-1}(\Delta)) \subset \Delta$.
   (The direction in which we want the maps $ h, g, \alpha$ and $\alpha'$
 to   combinatorially
 commute    is indicated by the first map  in the list. Thus saying
 that $\alpha',  h, g $ and $\alpha$ combinatorially commute we would mean that
   $(h\circ \alpha)(( \alpha'\circ  g)^{-1}(\Delta)) \subset \Delta$
   for every simplex $\Delta $ of $K'$.)
 Recall that a  map $h' : K \lo L'$ is
 a combinatorial lifting of $h$  to $L'$  if
 for every  simplex  $\Delta$ of $K'$   we have that
 $(\alpha' \circ h')( h^{-1}(\Delta)) \subset \Delta$.

 For   a simplicial complex $K$ and $a \in K$,   $st(a)$  denotes the union of all the simplexes
 of   $K$ containing  $a$.
  The following proposition    whose
 proof  is left to the reader      is
 a collection of   simple combinatorial  properties of    maps.

 \begin{proposition}
 %{\rm (\cite{l3})}
 \label{p3}${}$

   (i)     Let a compactum $X$  be  represented as the inverse limit
 $X  ={\rm  \invlim }K_i$ of  finite simplicial complexes $K_i$
 with  bonding maps   $h_{j}^i  :  K_{j} \lo K_i$.  Fix $i$ and   let
 $\omega : EW(K_i,k) \lo K_i$ be a resolution of $K_i$ which is suitable for $X$.
 Then there is
  a sufficiently large $j$  such that $h_j^i$ admits a combinatorial lifting to $EW(K_i,k)$.

  (ii)       Let  $h : K \lo K'$, $h' : K \lo L'$ and $\alpha' : L' \lo K'$ be maps of
  a simplicial complex $K'$   and  CW-complexes $K$ and $L'$  such that
  $h$ and $\alpha'$ are combinatorial and $h'$ is a combinatorial lifting of $h$.
 Then  there is a cellular approximation of $h'$ which is  also a combinatorial lifting of $h$.

 (iii)    Let $K$ and $K'$ be simplicial complexes,
 let    maps    $h : K \lo K'$, $g : L \lo L'$,  $\alpha : L \lo K$ and $\alpha' : L' \lo K'$
 combinatorially commute and let $h$ be  combinatorial.  Then

 $ g( \alpha^{-1}(st( x)) ) \subset
   \alpha'{}^{-1}(st( h (x)) )$ and
   $h(st(\alpha (z)))  \subset        st((\alpha' \circ g)(z))$  \\
   for every   $x \in K$ and  $z\in L$.

 \end{proposition}

  \end{section}

  \begin{section}{Proof of Theorem \ref{t1}}
 If $n\geq \dim_\z X$ then Theorem \ref{t1} follows from Theorem \ref{l3}
 (see also a remark at the end of this section).
Hence we may assume that   $n <  \dim_\z X$.  Then  for a finitely generated  abelian group
  $G$  the condition  $\dim_G X \leq n$ implies  that $G$ is torsion.  Thus we may
  assume  that  the groups   $G$ considered in the theorem   are  torsion  and   therefore
   the Bockstein basis $\sigma (G)$ of $G$ consists  only
  of groups of type $\z_p$ ($p$ is always assumed to be a prime number).

    Represent $X$ as   the inverse limit $X  = \invlim (K_i,h_i)$ of finite simplicial complexes $K_i$
 with combinatorial bonding maps   $h_{i+1}  :  K_{i+1} \lo K_i$   and the projections
 $p_i : X \lo K_i$ such that for every simplex $\Delta$ of $K_i$,  diam$(p_i^{-1}(\Delta)) \leq 1/i$.
  We will construct  by induction finite  simplicial complexes
  $L_i$   and maps      $g_{i+1}: L_{i+1} \lo L_i$,
 $\alpha_i : L_i \lo K_i$   such that    \\

 (a)  $L_i= K_i^{[n]}$   and $\alpha_i    : L_i \lo K_i$ is the inclusion.  The simplicial structure
 of $L_1$ is      induced from     $ K_1^{[n]}$         and
  the simplicial structure
 of $L_i$, $i >1$ is
   defined as   a  sufficiently small barycentric subdivision of   $ K_i^{[n]}$.
   We will refer to this simplicial structure while constructing standard resolutions
   of $L_i$.
    It is clear that   $\alpha_i$ is
 always a combinatorial map;

 (b)  the maps $h_{i+1}$, $g_{i+1}$,  $\alpha_{i+1}$  and $\alpha_i$   combinatorially commute.
 Recall that   this  means that
 for every simplex  $\Delta $ of $ K_{i}$,
 $(\alpha_i  \circ g_{i+1})((h_{i+1} \circ \alpha_{i+1})^{-1}(\Delta)) \subset \Delta$.  \\

  We will construct $L_i$ in such a way that   $Z=\invlim(L_i,g_i)$ will  admit
a map $r : Z \lo X$ such that  $Z$ and $r$ satisfy the conclusions of the theorem.
 Assume that the construction
 is completed for $i$.  We  proceed to $i+1$ as follows.

   Let
 $\dim_{\z_p}  X \leq k$, $2 \leq k  \leq n$ and let $f  :  N  \lo K(\z_p ,k)$ be a  cellular map  from
 a subcomplex $N$ of $L_i$, $  L_i^{[k]}  \subset  N$.
Let  $\omega_L: EW(L_i, k) \lo L_i^{}$
be the standard resolution  of $L_i$ for $f$.   We are  going to construct
from  $\omega_L: EW(L_i, k) \lo L_i^{}$ a resolution $\omega : EW(K_i,k)\lo K_i$ of
 $K_i$ suitable for $X$.  If $\dim K_i \leq  k$ set
 $\omega=\alpha_i \circ \omega_L      :  EW(K_i,k)=EW(L_i, k) \lo K_i$.

If $\dim K_i  >  k$ set
$\omega_{k}=\alpha_i \circ \omega_L :    EW_{k}(K_i,k)  =EW(L_i, k) \lo  K_i$
  and we will construct by induction  resolutions
 $\omega_j : EW_j(K_i, k) \lo K_i$, $k+1\leq j\leq \dim K_i$
  such that   $EW_{j}(K_i, k)$ is a subcomplex of
   $EW_{j+1}(K_i, k)$ and $\omega_{j+1}$ extends $\omega_j$  for every $k\leq j  < \dim K_i$.
 Note that saying that an $(n+1)$-cell is attached to a CW-complex by
    a map of degree  $p$ we mean  that the cell is attached  by a map $\phi$
    from the boundary $\s^n$ of the cell
  to the CW-complex  such that $\phi$    factors through   a map  $\s^n \lo \s^n$ of degree $p$.

    Assume that $\omega_j :   EW_j(K_i, k) \lo K_i$, $k\leq j <\dim K_i$
    is constructed.
    For every simplex $\Delta$ of $K_i$ of $\dim=j+1 $   consider  the subcomplex
    ${\omega_j}^{-1}(\Delta)$ of  $EW_j(K_i, k)$.  Enlarge
    ${\omega_j}^{-1}( \Delta)$    by
  attaching  cells  of $\dim =n+1  $     by maps of degree $p$
   in order  to kill     the elements of
   $p\pi_n({\omega_j}^{-1}( \Delta))=\{pa : a \in \pi_n({\omega_j}^{-1}( \Delta))\}$
   and attaching    cells
   of $ \dim > n+1  $ in order  to  get
    a subcomplex with trivial homotopy groups in
 $\dim   > n$.
 Let   $EW_{j+1}(K_i, k)$ be $EW_j(K_i, k) $ with all the cells  attached   for
 all  $(j+1)$-dimensional simplexes  $\Delta$ of $K_i$ and let
  $\omega_{j+1}:      EW_{j+1} (K_i, k)\lo K_i$ be an extension of
  $\omega_j$
      sending the interior points of the attached cells to the interior of       the corresponding
    $\Delta$.

      Finally denote  $EW(K_i, k)=EW_j(K_i, k)$ and $\omega=\omega_j : EW_j(K_i, k) \lo K_i$
    for $j=\dim K_i$.       Note that since we attach cells only of $\dim >n$,
    the $n$-skeleton of   $EW(K_i, k)$     coincides   with  the $n$-skeleton of $EW(L_i,k)$.

        Let us show that $EW(K_i, k)$ is suitable for $X$.
        Fix a simplex $\Delta$ of $K_i$ and denote $T= \alpha_i^{-1}(\Delta)$.
      First  note that  $T$  is $(n-1)$-connected,
       $\omega^{-1}(\Delta )$  is $(k-1)$-connected,
      $\pi_n(    \omega^{-1}(\Delta ))$ is $p$-torsion,
          $\pi_j(\omega^{-1}(\Delta ))=0$  for $j \geq  n+1$
    and $\pi_j(\omega^{-1}(\Delta ))=\pi_j (\omega_L^{-1}(T))$ for $ j \leq n-1$.

    By   Proposition \ref{p1},
        $\pi_j (\omega_L^{-1}(T))$
        is $p$-torsion for $j \leq n-1$.  Then
         $\pi_j(\omega^{-1}(\Delta ))$ is $p$-torsion for     $k \leq j \leq n$. Therefore
         by Bockstein Theory $\dim_{ \pi_j(\omega^{-1}(\Delta ))} X \leq \dim_{\z_p} X \leq k$,
         $k \leq j \leq n$
           and
        hence   by Proposition \ref{p2}, $\edim X \leq     \omega^{-1}(\Delta )$.   \\

          Thus we have shown that $EW(K_i, k)$ is suitable for $X$.
       Now   replacing  $K_{i+1}$ by
         $K_j$ with a sufficiently large $j$ we may   assume
    by (i)  of Proposition  \ref{p3}       that there is a combinatorial lifting
         of $h_{i+1}$ to $h'_{i+1} : K_{i+1} \lo  EW(K_i,k)$.
        By (ii) of  Proposition  \ref{p3} we  replace $h'_{i+1}$ by its  cellular approximation
         preserving
  the property of   $h'_{i+1}$ of being a      combinatorial lifting of $h_{i+1}$.

Then $h'_{i+1}$  sends the $n$-skeleton of $K_{i+1}$ to
the $n$-skeleton of $EW(K_i,k)$.  Recall that the  $n$-skeleton of $EW(K_i,k)$ is
contained in  $EW(L_i,k)$ and hence
 one can define
 $g_{i+1}=\omega_L \circ h'_{i+1}|_{K_{i+1}^{[n]}} : L_{i+1}= K_{i+1}^{[n]}\lo L_i$.
 Finally define a simplicial structure on $L_{i+1}$ to be a sufficiently small barycentric
 subdivision of      $K_{i+1}^{[n]}$  such that   \\

  (c)
 diam$g_{i+1}^j(\Delta) \leq 1/i$ for every simplex
 $\Delta$ in  $L_{i+1}$  and $j \leq i$   \\
     where
 $g^j_i=g_{j+1} \circ g_{j+2} \circ ...\circ g_i: L_i \lo L_j$.   \\

  It is easy to check that the properties (a) and (b)
 are satisfied.   \\

        Denote $Z=\invlim(L_i,g_i)$ and let $r_i : Z \lo L_i$ be the projections.

        Clearly $\dim Z \leq n$.
            For constructing
           $L_{i+1}$  we used an arbitrary map $f : N  \lo K(\z_p,k)$ such that
           $\dim_{\z_p} X \leq k$,  $ 2 \leq k \leq n$  and  $N$ is a subcomplex
 of $L_i$ containing $L_i^{[k]}$.
By  a standard reasoning described in detail   in the proof of Theorem 1.6,  \cite{l3}
one can show  that  choosing $\z_p$ and $f$  in an appropriate way
           for each $i$ we can achieve
that $\dim_{\z_p} Z \leq k$  for every integer $2\leq k \leq n$ and
 every $\z_p  $  such that      $\dim_{\z_p}  X \leq k $.
 Then by the Bockstein theorem
    $\dim_G Z \leq k$      for every    finitely generated  torsion abelian group $G$
   such that    $  \dim_G X \leq k $, $2 \leq k \leq n$.
       \\

   By (iii) of Proposition \ref{p3},        the property    (b) implies
    that for every   $x \in X$ and  $z\in Z$
   the following  holds:\\

 (d1)
  $ g_{i+1}( \alpha_{i+1}^{-1}(st( p_{i+1} (x) ) )) \subset
   \alpha_{i}^{-1}(st( p_i (x)) )$  and

   (d2)
   $h_{i+1}(st((\alpha_{i+1}\circ r_{i+1})(z)))  \subset        st((\alpha_{i}\circ r_{i})(z))$.
 \\

   Define a  map $r : Z \lo X$ by $r(z)=\cap \{ p_i^{-1}(    st((\alpha_{i}\circ r_{i})(z) ) ): i=1,2,... \}$.
  Then   (d2) implies
   that $r$ is indeed well-defined and continuous.

  The properties (d1) and (d2) also imply  that  for every $x \in X$  \\

 $r^{-1}(x)=\invlim ( \alpha_i ^{-1}(st(p_{i} (x))), g_i |_{\alpha_i ^{-1}(st(p_{i} (x)))})$\\
 where  the map  $ g_i |_{...} $ is considered as a map
      to $\alpha_{i-1} ^{-1}(st(p_{i-1} (x)))$. \\

Since $r^{-1}(x) $ is not empty for every $x \in X$,
 $r$ is a map onto. Fix $x \in X$ and  let us show that $r^{-1}(x)$  satisfies
the conclusions of the theorem.         First note    that     $M_i=st(p_{i} (x))$ is contractible .
Since    $T_i=\alpha_i ^{-1}(M_i)$ is  homeomorphic to the $n$-skeleton of     $st(p_{i} (x))$,
$T_i$ is  $(n-1)$-connected and hence
 $r^{-1}(x)$ is approximately $(n-1)$-connected
as the inverse   limit of $(n-1)$-connected finite simplicial complexes.

Let $\z_p$ be such that
 $\dim_{\z_p} X \leq k$,  $2\leq k \leq n$ and assume that $L_{i+1}$ is constructed
with help of   $f : N \lo K(\z_p, k)$.
Consider $T_{i+1}$ as the $n$-skeleton of $M_{i+1}$  and
 denote   \\

 $\beta$= the inclusion $:    T_{i+1} \lo   M_{i+1}$,

  $\tau=h'_{i+1}|_{...}       :   T_{i+1} \lo     \omega_L^{-1}(T_i)$,

 $\gamma =  h'_{i+1}|_{ ...} :     M_{i+1} \lo   \omega^{-1}(M_i)$  and

  $\kappa$=the inclusion$:    \omega^{-1}_L  (T_i) \lo     \omega^{-1}(M_i)$.\\

 Clearly $\gamma \circ  \beta = \kappa \circ \tau$.  Denote  by    $\beta^*,  \tau^*,\gamma^*$ and
 $\kappa^*$   the induced homomorphisms of the $n$-dimensional cohomology groups
modulo $\z_p$         of the corresponding spaces.
  Recall  that the $(n+1)$-cells of  $\omega^{-1}(M_i)$
not contained in $\omega_L^{-1}(T_i)$         are attached to  $\omega_L^{-1}(T_i)$       by maps
of degree $p$.
Then
$\kappa^*$ is an isomorphism.
 Since   $M_{i+1}$ is contractible,
 $\beta^*$ is the zero homomorphism.  Thus we obtain  that $\tau^*$ must also be
 the zero homomorphism.       Then, since  $g_{i+1}|_{...} : T_{i+1} \lo T_i$
 factors  through  $\tau$,  the map  $g_{i+1}|_{...}$ induces the     zero homomorphism  of  $H^n( T_{i}; \z_p)$
 and $H^n(T_{i+1}; \z_p)$.  Now  we may assume that  $\z_p$      appears in the construction
 for infinitely  many indices $i$.    Then $\check{H}^n(r^{-1}(x); \z_p)=0$
 and since  $  r^{-1}(x)$ is  approximately $(n-1)$-connected  and of $\dim \leq n$   we  get that
 $r^{-1}(x)$ is $\z_p$-acyclic.

Let $G$ be a finitely generated torsion abelian group  such that      $\dim_G X \leq n$.
 By the Bockstein theorem
 $\dim_{\z_p}  X \leq n$    for every   $\z_p \in \sigma(G)$.
 Then  $r^{-1}(x)$ is  $\z_p$-acyclic   for every  $p$ such that $\tor_p G \neq 0$.
 Hence   $r^{-1}(x)$ is  $G$-acyclic   and
 the theorem  follows.
\hfill $\Box$
     \\  \\
 {\bf Remark.}  It is easy to see that the proof of Theorem \ref{t1}  also works for
 $\z$  regarded as $\z_p$ with $p=0$.  This way one can avoid the use of Theorem \ref{l3}
  in the proof of Theorem \ref{t1}  and make  the proof  self-contained.

  \end{section}

Department of Mathematics\\
Ben Gurion University of the Negev\\
P.O.B. 653\\
Be'er Sheva 84105, ISRAEL  \\
e-mail: mlevine@math.bgu.ac.il\\\\
\end{document}